\documentstyle[12pt]{article}
\input{tcilatex}

\begin{document}

\vspace{1.cm}

\begin{center}
{\huge {\bf Problems on polygons and Bonnesen-type inequalities}} \footnote{%
Mathematics Subject Classification 51M10, 51M25, 52A40} \\[0pt]

\end{center}

\vspace{2.cm} {\large \ A. Raouf Chouikha}

\vspace{3.cm}

{\bf Abstract}\newline
{\it In this paper we are interested in some Bonnesen-type isoperimetric
inequalities for plane n-gons in relation with the two conjectures proposed
by P. Levy and X.M. Zhang.}\newline

\vspace{2.cm}

\section{Introduction}

As a well known result, for a simple closed curve ${\cal C}$ (in the
euclidian plane) of length $L$ enclosing a domain of area A, we have the
inequality

\begin{equation}  \label{1}
L^2 - 4\pi A \geq 0.
\end{equation}
Equality is attained if and only if this curve is a euclidean circle. This
means that among the set of domains of fixed area, the euclidean circle has
the smallest perimeter.\newline
The above inequality (\ref{1}) could be easily deduced from the Wirtinger
inequality

\begin{equation}
\int^{2\pi}_0|f^{\prime}(x)|^2 dx \geq \int^{2\pi}_0|f(x)|^2 dx,  \label{2}
\end{equation}
where $f(x)$ is a continuous periodic function of period $2\pi$ whose
derivative $f^{\prime}(x)$ is also continuous and $\int^{2\pi}_0 f(x) dx = 0$%
. Equality holds if and only if $f(x) = \alpha \cos x + \beta \sin x$ (see
Osserman's paper [O1]).\newline
For any curve ${\cal C}$ of length $L$ enclosing an area $A$, the quantity
\quad $L^2 - 4\pi A$\quad is called {\it the isoperimetric deficiency} of\ $%
{\cal C}$,\ because it decreases towards zero when \ ${\cal C}$ \ tends to a
circle. \newline
As an extension, Bonnesen proves [O1] that if  ${\cal C}$ is convex and
there exists a circular annulus containing ${\cal C}$ of thickness $d$, then we
have
$$4 \pi d^2 \leq L^2 - 4\pi A.$$
In fact, Fuglede showed that convexity is not a necessary hypothesis [F].%
\newline

There is a related isoperimetric inequality known as the Bonnesen :

\begin{equation}
L^2 - 4\pi A \geq \pi^2 (R - r)^2,  \label{3}
\end{equation}where\ $R$ \ is the circumradius and \ $r$ \ is the inradius of the curve \ $%
{\cal C}$. \newline
Note that if the right side of (\ref{3})  equals zero, then $R = r$.
This means that ${\cal C}$ is a circle and $L^2 - 4\pi A = 0$.\newline
More generally, inequalities of the form 
\begin{equation}
L^2 - 4\pi A \geq K,  \label{4}
\end{equation}
are called Bonnesen-type isoperimetric inequalities if  equality is only
attained for the euclidean circle. In the other words, $K$ is positive and
satisfies the condition

\[
K=0 \quad \mbox{ implies } \quad \ L^2 - 4\pi A = 0.
\]
(See [O2] for a general discussion and different generalisations).\newline

For an $n$-gon $\Pi_n$ (a polygon with $n$ sides) of perimeter $L_n$ and area $%
A_n$, the following inequality is known

\begin{equation}  \label{5}
L_n^2 - 4(n \tan{\frac{\pi}{n}}) A_n \geq 0.
\end{equation}
Equality is attained if and only if the $n$-gon is regular. Thus, if we
consider a smooth curve as a polygon with infinitely many sides, it appears
that inequality (\ref{1}) is a limiting case of (\ref{5}).\newline

Moreover, we know that an $n$-gon has a maximum area among all $n$%
-gons with the given set of sides if it is convex and inscribed in a circle.
Let \quad $a_1, a_2, ....,a_n$ \quad denote the lengths of the sides of \ $%
\Pi_n.$ \ For a triangle, Heron's formula gives the area

\[
A_3 = {\frac{1}{4}} L_3^2 \sqrt{ (1 - {\frac{2a_1}{{L_3}}})(1 - {\frac{2a_2}{%
{L_3}}})(1 -{\frac{2a_3}{{L_3}}})}.
\]
For a quadrilateral, the Brahmagupta formula gives a bound for the area

\[
A_4 \leq {\frac{1}{4}} L_4^2 \sqrt{ (1 - {\frac{2a_1}{{L_4}}})(1 - {\frac{%
2a_2}{{L_4}}})(1 - {\frac{2a_3}{{L_4}}})(1 - {\frac{2a_4}{{L_4}}})}
\]
with equality if and only if \ $\Pi_4$ \ can be inscribed in a circle.

\section{ Isoperimetric constants}

We can ask if it is possible to get an analogous formula for other plane  polygons (not %%@
necessarilly inscribed in a circle). More precisely, is
the area $A_n$ of the $n$-gon is close to the following expression ?

\begin{equation}
P_n = {\frac{ L_n^2}{{4}}} \sqrt{ (1 - {\frac{2a_1}{{L_n}}})(1 - {\frac{2a_2%
}{{L_n}}})(1 - {\frac{2a_3}{{L_n}}}).....(1 - {\frac{2a_n}{{L_n}}})}
\label{6}
\end{equation}

This question has been considered by many geometers who tried to
compare\quad $\displaystyle A_n$\quad with \quad $P_n$. One of them, P. Levy
[L] was interested in this problem and more precisely he expected the
following\newline

{\bf Conjecture 1 :} {\it Define the ratio ${\displaystyle \varphi_n = {\frac{%
A_n}{{P_n}}}}$. For any $n$-gon \ $\Pi_n$, with sides $a_1,a_2,....,a_n$
enclosing an area $A_n,$ and $P_n$ defined as above, this ratio verifies }

{\it 
\[
{\displaystyle a)\quad {\frac{e}{{\pi}}}\leq \varphi_n}\qquad \mbox{and} \qquad
b)\quad {\displaystyle \varphi_n\ \leq\ 1.}
\]
}\newline

{\bf Remark :} {\it Notice that part a) is obviously only valid for cyclic $n
$-gons, but part b) of Conjecture $1$ may be true for any $n$-gon. In
particular, for a triangle we have \ $\varphi_3 = 1$\ and for a
quadrilateral, \ $\varphi_4 \leq 1\ \ (\varphi_4 = 1$\ in the cyclic case).%
\newline
}

Conjecture 1 was originally motivated by study of cyclic $n$-gons. More
precisely, for regular $n$-gons we get\ ${\displaystyle a_i=L_n/n}$. The
associated value of $\varphi_n$ is given by

\begin{equation}
{\varphi_n}^0 = {\frac{1}{{{n \tan{\frac{\pi}{{n}}}}\ (1-{\frac{2}{{n}}}%
)^{n/2}}}} \label{7}
\end{equation}
and satisfies the inequalities of Conjecture 1. Moreover, we may verify that 
${\varphi^0}_n$ is a decreasing function in $n$.\newline

As we shall see below, the lower bound $\varphi_n \geq {\frac{e}{{\pi}}}$
seems to have more geometric interest than the upper one. Indeed, it allows
one to estimate the defect between any  $n$-gon $\Pi_n$ and the
regular one. This defect may be measured by the quotient

\begin{equation}
\tau_n = {\frac{{\varphi_n}}{{{\varphi_n}^0}}} \label{8}
\end{equation}
which  tends to $1$ whenever $\Pi_n$ is close to being regular.
Moreover, $\tau_n$ is related to a new Bonnesen-type inequality for plane
polygons.\newline

On the other hand, H.T. Ku, M.C. Ku and X.M. Zhang, ( [K.K.Z] and [Z]) have
been interested in this same problem. Their approach is quite different.
They consider the so called pseudo-perimeter of second kind $\hat{L}_n$
defined by

\begin{equation}
\hat{L}_n = L_n ({\frac{n}{{n-2}}}) [ (1 - {\frac{2a_1}{{L_n}}})(1 - {\frac{%
2a_2}{{L_n}}})(1 - {\frac{2a_3}{{L_n}}}).....(1 - {\frac{2a_n}{{L_n}}}%
)]^{1/n}.  \label{9}
\end{equation}
In fact, there is a relation between $\hat{L}_n$ and $P_n$ 
\begin{equation}
\hat{L}_n = ({\frac{n}{{n-2}}}) (4{P_n})^{\frac{2}{{n}}} {L_n}^{\frac{n-4}{{n%
}}}.  \label{10}
\end{equation}
X.M. Zhang ([Z] p. 196) has proposed the following\newline

{\bf Conjecture 2:} {\it For any cyclic $n$-gon \ $\Pi_n,$ \ we have }

{\it 
\[
{\hat{L}_n}^2 - 4 (n\tan{\frac{\pi}{{n}}}) A_n \geq 0.
\]
}
{\it Equality holds if and only if \ $\Pi_n$ \ is regular.}\newline

For any $n$-gon\quad $\Pi_n,$ \ we have the natural inequality \quad $\hat{L}%
_n \leq L_n.$\quad The equality \ $\hat{L}_n = L_n$ \ holds if and only if \ 
$\Pi_n$ \ is regular (see Lemma (4-6) of [Ch]).\newline

Moreover, it has been remarked by Zhang that Conjecture 2 implies Conjecture
(2-6) of [K,K,Z] concerning the 3-parameter family of pseudo-perimeters
denoted by \quad ${\cal L}_n[x,(n-1)y,{\frac{nz}{{2}}}]$\quad for any $n$%
-gon inscribed in a circle. They prove that \quad ${\hat L}_n {\leq} {\cal L}%
_n {\leq} L_n$\quad where \quad ${\cal L}_n[1,0,0] = L_n$ \quad and \quad $%
{\cal L}_n[0,0,1] = {\hat L}_n$. \newline

More generally, we also examine the following \newline

{\bf Problem 2':} {\it Let us consider a piecewise smooth closed curve $%
{\cal C}$ in the euclidean plane, of length \ $L$\ and area \ $A$.\ Let $%
(\Pi_n)_n$ be a sequence of $n$-gons approaching ${\cal C}$. ${{L}_n}$, ${%
\hat{L}_n}$ and $A_n$ are respectively the perimeter, the pseudo-perimeter and
the area of $\Pi_n$. Supposing that \quad ${\hat L} = \lim_{n\rightarrow
\infty} {\hat{L}_n}$ \quad exists, do we have the Bonnesen-type inequality}

\[
{\hat L}^2 - 4 \pi A \geq 0 \ ?
\]
Examples given below show that Problem 2' may have a solution.\newline

In this paper, we shall discuss these conjectures and  exhibit the special
role played by \quad $\tau_n = {\frac{{\varphi_n}}{{{\varphi_n}^0}}},$
\quad where \ $\varphi_n$\ and ${\varphi_n}^0$ are defined as above for any
cyclic $n$-gon, with sides\ $a_1,a_2,...,a_n$.\newline

Accordingly, we also introduce the ratio

\begin{equation}
\nu_n = ({\frac{L_n}{{\hat{L}_n }}})^{{\frac{n}{2}}-2}.  \label{11}
\end{equation}
We will describe some examples. As a consequence we propose a conjecture
which seems to be more appropriate than Conjecture 2. In particular, it
yields bounds for \ $\tau_n$. The theorem below shows that the position of $%
\tau_n$ \ compared with \ $1$ \ and \ $\nu_n$ \ gives partial answers to 
both the above conjectures.\newline
\smallskip

{\bf Theorem 1}:\newline
{\it Let \ $\tau_n = {\frac{{\varphi_n}}{{{\varphi_n}^0}}}$ \ and \ $\nu_n =
({\frac{L_n}{{\hat{L}_n }}})^{{\frac{n}{2}}-2}$ \ be the constants
associated to any cyclic n-gon \ $\Pi_n$ \ , with sides \ $a_1,a_2,...,a_n$%
. \quad $L_n$ \ and \ $\hat{L}_n$ \quad are respectively the perimeter and
the pseudo-perimeter. We then have \newline
}
{\it \noindent ({\rm i}) The inequality \quad $\tau_n \leq 1 $ \quad implies
conjecture $1$ b) and conjecture $2$. Moreover, this implication is strict.%
\newline
({\rm ii}) The inequalities \quad $1 \leq \tau_n \leq \nu_n$ \ imply
conjecture $1$ a) and conjecture $2$.\newline
({\rm iii}) The inequality \quad $\nu_n < \tau_n$ \quad contradicts
conjecture $2$.\newline
}
{\it In these three cases, Equality \ $1 = \tau_n = \nu_n$ \ holds if and
only if $\Pi_n$ is regular .}\newline

\noindent Case (i) of Theorem 1 implies in particular that

\[
{\frac{{\varphi_n}}{{{\varphi_n}^0}}} \leq 1 \leq ({\frac{L_n}{{\hat{L}_n }}}%
)^{{\frac{n}{2}}-2}.
\]
Case (ii) will be illustrated below by several examples. We hope that the
following hypothesis \quad $\varphi_n \leq {\varphi_n}^0$\quad will be
verified by an $n$-gon.\newline
As a corollary, we deduce from (ii) and (iii) that \quad $\tau_n \leq \nu_n$
\quad is equivalent to Conjecture 2.\newline

Consequently, we also obtain the following result.\newline

\newpage

{\bf Corollary 2}\newline
{\it Suppose \ $\tau_n \leq 1$ \ is verified by a cyclic n-gon; we then have
the following Bonnesen-type isoperimetric inequality : }

{\it 
\[
\hat{L}_n^2 - 4(n \tan{\frac{\pi}{n}}) A_n \ \geq \ \hat{L}_n^2 (1 - {\frac{{%
\varphi_n}}{{{\varphi_n}^0}}}).
\]
}
{\it Equality holds if and only if $\Pi_n$ is regular (i.e. \quad ${\varphi_n%
} = {\varphi_n}^0$).\newline
Moreover, this inequality implies Conjecture $2$\ .}\newline

\section{ Proofs}

{\bf 1.}\quad Let \quad $L_n,\ {\hat L}_n,\quad A_n$\quad be respectively
the perimeter, pseudo-perimeter and area of any polygon\quad $\Pi_n$\quad
as defined in the preceding section. The sides are of lengths \quad $%
a_1,a_2,...,a_n$. Consider ratio ${\displaystyle \varphi_n = {\frac{A_n}{{%
P_n}}}}$, \quad where

\begin{equation}
P_n = {\frac{ L_n^2}{{4}}} \sqrt{ (1 - {\frac{2a_1}{{L_n}}})(1 - {\frac{2a_2%
}{{L_n}}}) .....(1 - {\frac{2a_n}{{L_n}}})}= {\frac{1}{{4}}}({\frac{n-2}{n}}%
)^{\frac{n}{2}}({\hat L_n})^{\frac{n}{2}}(L_n)^{\frac{4-n}{2}}.  \label{12}
\end{equation}
Then we get expression

\[
\displaystyle \varphi_n = ({\frac{n-2}{n}})^{-{\frac{n}{2}}}\ {\frac{4 A_n}{{%
L_n^2}}}\ ({\frac{L_n}{{\hat L_n}}})^{\frac{n}{2}}.
\]

After simplification, we have

\[
{\frac{{\varphi_n}}{{{\varphi_n}^0}}} = {\frac{4(n \tan{\frac{\pi}{n}}) A_n}{%
{L_n^2}}}\ ({\frac{L_n}{{\hat L_n}}})^{\frac{n}{2}}.
\]
Consequently, we obtain a relation between $\tau_n$ and $\nu_n:$ 
\begin{equation}
\tau_n = {\frac{4(n \tan{\frac{\pi}{n}}) A_n}{{L_n^2}}}\ (\nu_n)^{\frac{n}{{%
n-4}}}  \label{13}
\end{equation}

or

\begin{equation}
\tau_n = {\frac{4(n \tan{\frac{\pi}{n}}) A_n}{{\hat L_n^2}}}\ \nu_n.
\label{14}
\end{equation}
This proves that Conjecture 2 is equivalent to the inequality

\[
\tau_n \leq \nu_n.
\]
Furthermore, since examples given below verify condition (ii) of Theorem 1, $%
\ 1 \leq \tau_n \leq \nu_n$, we may deduce that the implication (i) is
necessarily strict.\newline

Moreover, $\displaystyle {\varphi_n} \leq {\varphi_n}^0 \leq 1$ implies that 
$\displaystyle ({\frac{L_n}{{\hat L_n}}})^{\frac{n}{2}} \leq {\frac{{L_n}^2}{%
{4(n \tan{\frac{\pi}{n}})A_n}}}.$ The latter implies

\[
\nu_n \leq {\frac{{{\hat L}_n}^2}{{4(n \tan{\frac{\pi}{n}})A_n}}},
\]
which is equivalent to Conjecture 2, since $\nu_n \geq 1 .$\newline

We may deduce from the above some necessary conditions satisfied by $\tau_n$%
. Indeed, from (\ref{13}) and (\ref{14}), the ratio should verify the
inequalities

\[
\tau_n \leq {\nu_n}^{\frac{n}{{n-4}}}, \qquad \tau_n \geq {\frac{4(n \tan{%
\frac{\pi}{n}}) A_n}{{\hat L_n^2}}} \geq {\frac{4(n \tan{\frac{\pi}{n}}) A_n%
}{{\ L_n^2}}}.
\]

All the equalities are attained only if\ $\nu_n = 1,$ \ which corresponds
to the regular polygon. Theorem 1 is proved.\newline

{\bf 2.}\quad We prove now Corollary 2. Since \ $\nu_n \geq 1$\ we may
deduce from (14) the following :

\[
{\frac{{\varphi_n}}{{{\varphi_n}^0}}} \geq {\frac{4(n \tan{\frac{\pi}{n}})
A_n}{{\ \hat{L}_n^2}}}.
\]
We then get

\[
1 - {\frac{4(n \tan{\frac{\pi}{n}}) A_n}{{\ \hat{L}_n^2}}} \geq 1 - {\frac{{%
\varphi_n}}{{{\varphi_n}^0}}} \geq 0.
\]

\noindent Thus, 

\[
0 \leq {\hat L}_n^2(1 - \tau_n) \leq {\hat L}_n^2 - 4(n\tan{\frac{\pi}{n}})
A_n.
\]

So, we have proved the first part of Corollary 2. This inequality implies
obviously ${\hat L}_n^2 - 4n \tan{\frac{\pi}{n}} A_n \geq 0$, i.e.
Conjecture 2. Moreover, it is clear that equality is attained for the regular polygon $%
\Pi_n$.\newline

Conversely, suppose we have

\[
{\hat L}_n^2 - 4n \tan{\frac{\pi}{n}} A_n \ = \ {\hat L}_n^2 (1 -{\frac{{%
\varphi_n}}{{{\varphi_n}^0}}}).
\]
This is equivalent to

\begin{equation}
\tau_n = {\frac{4n \tan{\frac{\pi}{n}} A_n}{{\ {\ L}_n^2}}}.  \label{15}
\end{equation}
That means \quad $\nu_n = 1,$\quad i.e. \ $\Pi_n$ \ is regular (see [Ch],
Lemma(4.6) ).\newline

\noindent {\bf Remark}\newline
Under the hypothesis of Corollary 2, suppose in addition, that ${\tau_n}{%
\nu_n} \geq 1$. We then obtain a better Bonnesen-type isoperimetric
inequality

\[
{\hat L}_n^2(1 - \tau_n) \leq {\hat L}_n^2(1 - {\frac{1}{{\nu_n}}}) \leq {%
\hat L}_n^2 - 4n \tan{\frac{\pi}{n}} A_n.
\]

\section{ Some special polygons}

In this part, we shall see that Hypothesis (ii) of Theorem 1, which implies
Conjecture 2, is in fact verified by many examples.

\subsection{Example 1}

It is true in particular for the Macnab polygon, which is a cyclic equiangular alternate-sided %%@
2$n$-gon
with $n$ sides of length $a$ and $n$ sides of length $b$.  This polygon was originally used as an %%@
example by [K,K,Z] and by [Z], to
test their conjectures.\newline
In fact, we can do better by the following result:\newline

{\bf Proposition 1}\newline
{\it Let $\Pi _{n,n}$ be a cyclic 2$n$-gon with $n$ sides of length $a$
alternatively with $n$ sides of length $b$ and $\varphi _{n,n}$ its
associated function. Then, we have}

\[
1\leq \frac{\varphi _{n,n}}{\varphi ^{0}}\leq (\frac{L_{n,n}}{\hat{L}_{n,n}}%
)^{\frac{n}{2}-2}.
\]

{\bf Proof} 

 A direct calculation gives the expression

\[
\varphi _{n,n}=\frac{[(a^{2}+b^{2})\cos \frac{\pi }{n}+2ab]}{n\sin \frac{%
\pi }{n}(a+b)^{2}[1-\frac{2}{n}+ \frac{4ab}{n^{2}(a+b)^{2}}]^{\frac{n}{2}}}.
\]
This follows from expressions for $A_{n,n}$\quad and \quad $\hat{L}_{n,n}$\quad
calculated by [Z].\\
Indeed, one gets

\begin{eqnarray*}
A_{n,n}=\frac{n}{4\sin \frac{\pi }{n}}[(a^{2}+b^{2})\cos \frac{\pi }{n}
+2ab], \\
\hat{L}_{n,n}^{2}=\frac{n^{2}}{(n-1)^{2}}[%
n(n-2)(a^{2}+b^{2})+n^{2}ab+(n-2)^{2}ab]. \\
\end{eqnarray*}

Furthermore, it is easy to see that $\varphi _{n,n}$ may be written

\begin{eqnarray*}
\varphi _{n,n}=\frac{[\cos \frac{\pi }{n}+\frac{2ab}{(a+b)^{2}}(1-\cos 
\frac{\pi }{n}]}{n\sin \frac{\pi }{n}[1-\frac{2}{n}+\frac{4ab}{%
n^{2}(a+b)^{2}}]^{\frac{n}{2}}}.
\end{eqnarray*}
Thus, 
\begin{eqnarray*}
\frac{\varphi _{n,n}}{\varphi ^{0}}=\frac{(1-\frac{1}{n})^{n}[\cos \frac{%
\pi }{n}+\frac{2ab}{(a+b)^{2}}(1-\cos \frac{\pi }{n}]}{ 1-\frac{2}{n}+%
\frac{4ab}{n^{2}(a+b)^{2}}]^{\frac{n}{2}}(1+\cos \frac{\pi }{n})},
\end{eqnarray*}
which can be expressed as follows :

\begin{eqnarray*}
\frac{\varphi _{n,n}}{\varphi ^{0}}=(1+\frac{E-1}{n(n-1)})^{-\frac{n}{2}%
}[1+(E-1) \frac{1-\cos \frac{\pi }{n}}{1+\cos \frac{\pi }{n}}], \mbox{
where } E=\frac{4ab}{(a+b)^{2}}.
\end{eqnarray*}
We can see easily that $\frac{\varphi _{n,n}}{\varphi ^{0}}\geq 1$, since $%
0\leq E\leq 1$. \newline

\subsection{Example 2}

Let  $\Pi _{n}^{0}$ denote the regular $n$-gon whose sides $a_{i\quad }^{0}$
are subtended by angles $\frac{\pi}{n}; i=1,...,n$. Consider a polygon $\Pi
_{n}^{\varepsilon }$ obtained from $\Pi _{n}^{0}$ by variations of $%
a_{1},a_{2}$ which are subtended respectively by $\frac{\pi }{n}-\varepsilon$
and $\frac{\pi}{n}+\varepsilon$.\ The other sides of length $a_{i}^{0}
(3\leq i\leq n)$ are unchanged. We prove that hypothesis (ii) is verified
by \ $\Pi _{n}^{\varepsilon }$.\newline

{\bf Proposition 2}

{\it Let $\Pi _{n}^{\varepsilon }$ be the $n$-gon defined above for $n\geq 4$,
 $\varphi _{n}^{\varepsilon }$ being its associated function. Then, for $%
\varepsilon > 0 $ small, we have $1\leq \frac{\varphi _{n}^{\varepsilon }}{%
\varphi _{n}^{0}}\leq (\frac{L_{n,n}}{\hat{L}_{n,n}})^{\frac{n}{2}-2}.$}%
\newline

Thus, it seems that the function $\varphi(a_{1,}a_{2,}....,a_{n})$
for an $n$-gon possesses a local minimum for the regular polygons.%
\newline

{\bf Proof}

Let $L_{n}^{\varepsilon },\hat{L}_{n}^{\varepsilon},A_{n}^{\varepsilon }$ be
respectively perimeter, pseudo-perimeter and enclosing area of the polygon $%
\Pi _{n}^{\varepsilon } $ defined above. We get $a_{1}=2R\sin (\frac{\pi }{n}%
-\varepsilon )$ and $a_{2}=2R\sin (\frac{\pi}{n} +\varepsilon )$. After
calculation, we obtain the following expression

\begin{eqnarray*}
A_{n}^{\varepsilon }=A_{n}^{0}(1-\frac{4\varepsilon ^{2}}{n})\quad \mbox{and}\quad
L_{n}^{\varepsilon }=L_{n}^{0}(1-\frac{ \varepsilon ^{2}}{n}).
\end{eqnarray*}
On the other hand,
\begin{eqnarray*}
(1-\frac{2a_{1}}{L_{n}^{\varepsilon }})(1-\frac{2a_{2}}{L_{n}^{\varepsilon }}%
)&=&1-\frac{4}{n}[ 1-\varepsilon ^{2}(\frac{1}{2}-\frac{1}{n})]+\frac{%
4R^{2}[\sin ^{2} \frac{ \pi }{n}-\varepsilon ^{2})}{(L_{n}^{0})^{2}(1-\frac{%
\varepsilon ^{2}}{n})} \\
&=&(1-\frac{2}{n})^{2}[1+\frac{4\varepsilon ^{2}}{(n-2)^{2}}(\frac{n-2}{2}+ 
\frac{2}{n}-\frac{1}{\sin ^{2} \frac{\pi }{n}})].
\end{eqnarray*}
Also, we get $(1-\frac{2a_{i}^{0}}{L_{n}^{\varepsilon }})=(1-\frac{2}{n})(1-%
\frac{2\varepsilon ^{2}}{n(n-2)})$.\\
After simplification, we find the expression

\[
\frac{\varphi _{n}^{\varepsilon }}{\varphi _{n}^{0}}=1-\frac{2\varepsilon
^{2}}{(n-2)^{2}}[\frac{(n-2)^{2}}{2n}+\frac{n-2}{2}+\frac{2}{n}-\frac{1}{%
\sin ^{2}\frac{\pi }{n}}],
\]
which verifies $\frac{\varphi_{n}^{\varepsilon }}{\varphi _{n}^{0}}\geq 1.$%
\newline
Notice that the factor  $\varepsilon^2$ vanishes for $n=4$.\newline

From the expression

\[
(\frac{L_{n,n}}{\hat{L}_{n,n}})^{\frac{n}{2}-2} = 1 - {\frac{2(n-4)}{{%
n(n-2)^2}}}\varepsilon^2 [-\frac{(n-2)^{2}}{2n}+\frac{n-2}{2}+\frac{2}{n}-%
\frac{1}{\sin ^{2}\frac{\pi }{n}}],
\]

we also prove that $\frac{\varphi_{n}^{\varepsilon }}{\varphi _{n}^{0}}\leq (%
\frac{L_{n,n}}{\hat{L}_{n,n}} )^{\frac{n}{2}-2}.$\newline

\section{ Levy's polygons}

In this part, we discuss the connexion between Conjecture 1 and some
Bonnesen-type inequalities by using examples. Some $n$-gons satisfy Conjecture 1
 without being regular. P. Levy has remarked
on particular properties of the function $\varphi_{n}$ which depends 
 on the lengths of the sides

\[
\varphi _{n}=\varphi _{n}(a_{1},a_{2},a_{3},....,a_{n}).
\]
Indeed, he noticed that $\varphi _{n}$ is a bounded algebraic symmetric
function. Its bounds does not depend on $n$ and it should verify the equality

\begin{eqnarray*}
\varphi _{n}(a_{1},a_{2},a_{3},....,0)=\varphi
_{n-1}(a_{1},a_{2},a_{3},....,a_{n-1}).
\end{eqnarray*}
Consequently, we deduce that

\[
\tau_{n-1}(a_1,a_2,....,a_{n-1}) < \tau_n(a_1,a_2,....,a_{n-1},0).
\]

{\bf 5.1}\quad Also, P. Levy tried to find these bounds and tested Conjecture 1 on a
special curve polygon denoted by $\Pi (\alpha )$, inscribed in the euclidean
circle of radius 1. It is bounded by a circular arc with length $2(\pi -
\alpha)$, and a chord of length $l=2\sin \alpha$, where $0\leq \alpha \leq
\pi$.\\
$\Pi (\alpha )$ can be considered as limit of an $(n+1)$-gon with $n$ sides
of length $2\sin \frac{2\pi }{n}$ while only one has a fixed length $%
l=2\sin \alpha $. Let $\varphi _{n}(\alpha )$ be the corresponding ratio and 
$\varphi (\alpha )$ its limit value when $n$ tends to infinity. In this
case, ${\varphi^0}_\infty={\frac{e}{\pi }}$ is the limit value of ${\varphi^0%
}_n = {\frac{1}{{n \tan{\frac{\pi}{{n}}}\ (1-{\frac{2}{{n}}})^{n/2}}}}$.%
\newline

We get the following\newline

{\bf Proposition 3}\newline
{\it Let $L(\alpha ),\ {\hat L}\alpha ),\ A(\alpha )$ be respectively the
perimeter, the pseudo-perimeter and the enclosing area of the ``polygon'' $%
\Pi (\alpha )$, with $0\leq \alpha \leq \pi $. We then obtain the
inequalities\newline
}
{\it \noindent a) $1\leq \frac{\varphi (\alpha )}{\varphi ^{0}}\leq \frac{%
\pi }{e}\sqrt{\frac{e}{3}}$ with $\varphi (\pi )=\frac{e}{\pi }$  and $%
\varphi (0)=\sqrt{\frac{e}{3}}$.\newline
b) $\hat{L}^{2}(\alpha )-4\pi A(\alpha )\geq 0$.\newline
}
{\it Equality holds if and only if $\alpha = 0$.}\newline

\smallskip

Thus, we may deduce that $\Pi (\alpha )$ verifies Conjecture 1 and Problem
2'.\newline

{\bf Proof }

We may calculate the exact value of the function $\varphi (\alpha)$. We
refer for that to P. Levy's papers [L] and [Ch] for details. Here $%
L=2(\alpha +\sin \alpha )$ and $A=\alpha -\sin \alpha \cos \alpha$ \quad so that

\begin{eqnarray*}
\varphi (\alpha )=\frac{(\pi-\alpha -\sin \alpha \cos \alpha)}{(\pi-\alpha
+\sin \alpha )^{\frac{3}{2}}\sqrt{\pi-\alpha -\sin \alpha }}e^{\frac{%
\pi-\alpha}{\pi-\alpha+\sin \alpha }}.
\end{eqnarray*}

\[
\tau = \frac{\pi (\pi-\alpha -\sin \alpha  \cos \alpha)e^{\frac{\sin \alpha }{
\pi-\alpha+\sin \alpha }}}{ (\pi-\alpha +\sin
\alpha )^{\frac{3}{2}}\sqrt{\pi-\alpha -\sin \alpha }},
\]

and

\[
\nu = \sqrt{{\frac{(\pi-\alpha +\frac{1}{2}\sin \alpha ) }{%
(\pi-\alpha +\sin \alpha )}}} e^{\frac{\pi-2\alpha-\sin \alpha }{\alpha +\sin \alpha }}.
\]

Thus, for $0\leq \alpha \leq \pi$ we obtain the double inequality ( [Ch],
Proposition(2.1) )

\begin{eqnarray*}
\frac{e}{\pi }\leq \varphi (\alpha )\leq \sqrt{\frac{e}{3}}.
\end{eqnarray*}

These inequalities may also be verified by {\it Mathematica}. On the other
hand, we may also deduce the expression $\frac{4\pi A}{\hat{L}^{2}}$ in
terms of $\alpha$ :

\begin{eqnarray*}
\frac{4 A}{\hat{L}^{2}}= \frac{(\pi - \alpha -\frac{1}{2}\sin 2\alpha )%
}{(\pi - \alpha +\sin \alpha )^{2}}.
\end{eqnarray*}

We can prove easily that the right side of the above expression is a
decreasing function of $\alpha$, and for $\alpha = 0$, its value is $1$. We
then obtain part b) of Proposition 3.\newline

{\bf 5.2}\quad P. Levy considered also another curvilinear polygon. Denote by $\Pi(\alpha
,\theta )$ the polygon obtained from $\Pi (\alpha)$ by replacing the side with
length $l=2\sin \alpha$ by two sides. One of them has a length $2\sin \theta$%
. Then we get the expression of the perimeter and the area of the new
polygon $\Pi (\alpha ,\theta )$ :

\begin{eqnarray*}
L(\alpha ,\theta )=2[\pi - \alpha +\sin\theta +\sin (\alpha +\theta )], \\
A(\pi - \alpha ,\theta )=\pi - \alpha+\sin \alpha\cos (\alpha +2\theta ), \\
 0\leq \theta \leq \alpha.
\end{eqnarray*}
For $\theta =0$ we get of course, $\Pi (\alpha ,0)\equiv \Pi (\alpha ).$%
\newline

{\bf Proposition 4}\newline
{\it Let $L(\alpha ,\theta ),\ \hat{L}(\alpha ,\theta ),\ A(\alpha ,\theta)$
be respectively the perimeter, the pseudo-perimeter and the enclosing area
of the ``polygon'' $\Pi (\alpha ,\theta )$, with $0\leq \alpha \leq \pi $,
and $0\leq \theta \leq \pi-\alpha .$ We then obtain the inequalities\newline
}
{\it \noindent a) $\varphi (\alpha ,\theta _{0})\leq \varphi (\alpha ,\theta
)\leq \varphi (\alpha ,\frac{\pi -\alpha }{2})\leq 1$ for certain $\theta
_{0}>0$.\newline
b) $1\leq \frac{\varphi (\alpha ,\frac{\pi -\alpha }{2})}{\varphi ^{0}}\leq 
\frac{\pi }{e}$ with $\varphi (\pi ,0)=1$  and $\varphi (0,{\frac{\pi}{{2}}}%
)={\frac{\pi}{e}}$.\newline
c) $\hat{L}^{2}(\alpha ,\frac{\pi -\alpha }{2})-4\pi A(\alpha ,\frac{%
\pi-\alpha }{2})\geq 0$.\newline
}
{\it Equality holds if and only if $\alpha =\pi$.}\newline

{\bf Proof}\newline
We calculate the following expression for the function $\varphi (\alpha
,\theta)$ defined above 
\begin{eqnarray*}
\varphi (\alpha ,\theta )=\frac{\alpha -\sin \alpha \cos (\alpha +2\theta ) 
}{[\alpha +\sin \theta +\sin (\alpha +\theta )]\sqrt{\alpha ^{2}-4\sin ^{2} 
\frac{\alpha }{2}\cos ^{2}(\frac{\alpha }{2}+\theta )}}e^{\frac{\alpha }{%
[\alpha +\sin \theta +\sin (\alpha +\theta )]}}
\end{eqnarray*}
The details are given in [L] and [Ch]. In particular, for $0\leq\alpha \leq
\pi$ we have seen ([Ch], proposition (3-1)) that $\varphi (\alpha ,\theta)$
admits a maximum $\theta _{0}=\frac{\pi -\alpha }{2}$ and two minima $\theta
_{1},\theta _{2}$ symmetric with respect to $\theta _{0}$, such that $%
\varphi (\alpha ,\theta _{1})=\varphi (\alpha ,\theta _{2}).$\newline
Moreover, we may prove that $\frac{\varphi(\alpha ,\frac{\pi -\alpha }{2})}{%
\varphi ^{0}}$ is a decreasing function, \newline $\varphi (\pi ,0)=1$,
and $\varphi (0,\frac{\pi }{2})=\frac{\pi }{e}$.\newline

Furthermore, after simplifying the expression $\frac{4\pi A(\alpha ,\frac{%
\pi -\alpha }{2})}{\hat{L}^{2}(\alpha ,\frac{\pi -\alpha }{2})}$ we find the
following :

\begin{eqnarray*}
4 {\frac{A(\alpha ,\frac{\pi -\alpha }{2})}{{{\hat{L}^{2}(\alpha ,\frac{%
\pi -\alpha }{2})}}}}=  \frac{\pi - \alpha +\sin \alpha}{(\pi - \alpha
+2\cos \frac{\alpha }{2})^{2}}.
\end{eqnarray*}
We may verify that a such function is decreasing and is less than $1.$ We
have thus proved part c) of Proposition 4.\newline

{\bf Remark :} There are two possibilities for the ``polygon'' $\Pi (\alpha
) $ (considered as a limit of an $(n+1)$-gon) with only one side of length $%
l=2\sin \alpha$. The center of the circumscribed circle is inside the domain
bounded by the $(n+1)$-gons. In this case, $\sum 2\theta _{i}=2\pi$, where $%
2\theta _{i}$ is the subtended angle of the side $a_{i}$. The second case is
arises by transposing \ $\alpha$ \ with \ $\pi - \alpha$.\ So,  the center of
the circomscribed circle is outside the domain bounded by the $(n+1)$-gons,
and we have $\sum 2\theta_{i}<2\pi .$ This fact have been underlined by
P.Levy. Of course, it is only in the first case, that the isoperimetric inequality is
optimal.

\section{Concluding remark}

Thus, it is natural to expect that the hypothesis (ii) of Theorem 1 is
verified for any cyclic n-gon. We then may propose the following\newline

{\bf Conjecture 3} : {\it For any n-gon $\Pi_n,$ we have the inequalities }

{\it 
\[
1 \leq \tau_n \leq \nu_n,
\]
}

{\it with $1 = \tau_n = \nu_n$ if and only if $\Pi_n,$ is regular.}\newline

Obviously, this implies Conjecture 2 and Conjecture 1 a). Thus, conjecture 3
appears to be more significant than the previous conjectures. Notice that by
Theorem 1,

\[
\nu_n = 1 \Rightarrow \tau_n = 1.
\]
To investigate in this direction, we can see for example expression (13). We
may deduce anyway that ${\displaystyle  {\frac{{\nu_n}^{\frac{n}{{n-4}}%
}}{{\tau_n}}}}\geq 1$, which is equivalent to

\[
{\frac{{\varphi_n}}{{{\varphi_n}^0}}}({\frac{{\hat L}_n}{{L_n}}})^{\frac{n}{2%
}} \leq 1 \mbox{ and } {\frac{1}{{\nu_n^{\frac{4}{{n-4}}}}}} \leq {\frac{%
\nu_n}{{\tau_n}}}.
\]
Equality holds if and only if the $n$-gon is  regular. This gives an
upper bound for $\tau_n$.\newline
Actually, by using the Bonnesen-style inequalities of X.M. Zhang [Z], we can
improve it. More precisely, we have

\[
{\frac{{\tau_n}}{{{\nu_n}^{\frac{2n}{{n-4}}}}}} \leq 1 - ({\frac{2R n \sin{%
\frac{\pi}{n}}}{{L_n}}} -1)^2. 
\]
Also,

\[
{\frac{{\tau_n}}{{{\nu_n}^{\frac{2n}{{n-4}}}}}} \leq 1 - (1 - {\frac{2 r n
\tan{\frac{\pi}{n}}}{{L_n}}})^2.
\]
Here $R$ \ and \ $r$\ are the circumradius and inradius, respectively.%
\newline
Moreover, we get the following lower bound :

\[
({\frac{2 r n \tan{\frac{\pi}{n}}}{{L_n}}})^2 \leq {\frac{{\tau_n}}{{{\nu_n}%
^{\frac{2n}{{n-4}}}}}},
\]
which implies in particular, that

\[
{\nu_n}^{\frac{n+4}{{n-4}}}({\frac{2 r n \tan{\frac{\pi}{n}}}{{L_n}}})^2
\leq {\frac{{\tau_n}}{{\nu_n}}}.
\]

\newpage
{\bf Acknowledgment :} \ I would like to thank Professor Coxeter for his interest in my paper and %%@
his corrections.
\vspace{1cm}

{\bf REFERENCES}\newline

[Ch] R. Chouikha {\it Probleme de P. Levy sur les polygones articules}%
\newline
C. R. Math. Report, Acad of Sc. of Canada, vol 10, p. 175-180, 1988. \newline

[F] B. Fuglede {\it Bonnesen inequality for the isoperimetric deficiency of
closed curves in the plane} \qquad Geometriae Dedicata, vol 38, p. 283-300,
1991.\newline

[K.K.Z] H.T. Ku, M.C. Ku, X.M. Zhang {\it Analytic and geometric
isoperimetric inequalities} \qquad J. of Geometry, vol 53, p.100-121, 1995.%
\newline

[L] P. Levy {\it Le probleme des isoperimetres et des polygones articules}%
\newline
Bull. Sc. Math., 2eme serie, 90, p.103-112, 1966.\newline

[O1] R. Osserman {\it The isoperimetric inequalities}\newline
Bull. Amer. Math. Soc., vol 84, p.1182-1238, 1978.\newline

[O2] R. Osserman {\it Bonnesen-style isoperimetric inequalities}\newline
Amer. Math. Monthly, vol 1, p. 1-29, 1979.\newline

[Z] X.M. Zhang {\it Bonnesen-style inequalities and pseudo-perimeters for
polygons}\qquad J. of Geometry, vol 60, p.188-201, 1997.\newline

\vspace{2cm}

\begin{center}
University of Paris-Nord\\[0pt]
LAGA,CNRS UMR 7539\\[0pt]
Av. J.B.Clement\\[0pt]
Villetaneuse F-93430\\[0pt]
{\small e-mail: chouikha@math.univ-paris13.fr}
\end{center}

\end{document}